\documentclass[10pt,notitlepage,twoside]{article}
\usepackage[latin1]{inputenc}
\usepackage[T1]{fontenc}
\usepackage[english]{babel}
\usepackage{latexsym}
\usepackage{amsmath}
\usepackage{amsfonts}
\usepackage{amssymb}

\newcommand{\be}{\begin{equation}}
\newcommand{\ee}{\end{equation}}
\newenvironment{pf}{\noindent{\bf Proof.}\enspace}{
\hfill $\square$ \medskip}
\newcommand{\R}{\mathbb{R}}

\newcommand{\N}{\mathbb{N}}

\newtheorem{thm}{Theorem}[section]

\newtheorem{lem}{Lemma}[section]
\newtheorem{rem}{Remark}[section]
\newtheorem{cor}{Corollary}[section]
\newtheorem{df}{Definition}[section]

\DeclareFontFamily{T1}{cmr}{\hyphenchar\font=-1}

\numberwithin{equation}{section}
\DeclareFontFamily{T1}{cmr}{\hyphenchar\font=-1}
\usepackage{color}
\usepackage{latexsym}

\begin{document}
\title{On existence theorems of a  functional  differential
equations in  partially ordered Banach algebras }
\date{}
\maketitle

\par\vskip0.3cm


\begin{center} \emph{Amor Fahem$^{(1)}$, Aref Jeribi $^{(1)}$ and  Najib Kaddachi $^{(2)}$  }
\end{center}\vspace{3mm}


\begin{center}
\emph{$^{(1)}$ Department of Mathematics. University of Sfax.
Faculty of Sciences of Sfax. \\ Soukra Road Km 3.5 B.P. 1171, 3000, Sfax,
Tunisia}\\
E-mail: amorfahem.edp@gmail.com\\
E-mail: Aref.Jeribi@fss.rnu.tn

\end{center}

\begin{center}
\emph{$^{(2)}$ University of Kairouan. Faculty of Science and Technology of Sidi Bouzid. Agricultural University
City Campus $-$ 9100, Sidi Bouzid, Tunisia}\\
E-mail: najibkadachi@gmail.com

\end{center}

%
\vskip0.2cm
\begin{center}
 \vspace{3mm} \hspace{.05in}\parbox{5.5 in}{{\bf\small Abstract.}
  In this paper we are concerned with existence results for a coupled system of quadratic functional differential equations.  This system is reduced to a fixed point problem for a $2\times 2$ block
operator matrix with nonlinear inputs. To prove the existence we are established some fixed point theorem of Dhage's type for the block matrix operator acting in partially ordered Banach algebras.



  {\small \sloppy{}}}
\end{center}

\begin{center}
 \vspace{3mm} \hspace{.05in}\parbox{5.5 in}{{\bf\small Keywords:}
 {\small \sloppy{Partially ordered Banach algebra, Fixed point theory, Partial measure of noncompactness, differential equations, block Operator matrix.}}}
\end{center}
 \bigskip
\vspace{3mm} \hspace{.05in}\parbox{5.5 in}{
}

\section{Introduction}
\hspace*{20pt}The theory of fixed point is one of the most powerful and most fruitful tools of modern
mathematics and can consider a fundamental material of non-linear analysis. In recent years a number of
excellent monographs and surveys by distinguished authors about fixed point theory have appeared such as,
\cite{agrawel1, agrawel2, browder,caballero,dhage5, guolak}. Based on the fact that the Banach spaces are the fundamental underlying spaces on linear and nonlinear analysis, it
leads us to consider the following problem: if a Banach space is equipped with an ordering structure, partial
order or lattice, this Banach space becomes a partially ordered Banach space. Then, when we solve some problems
on this Banach space, in addition to the topological structure and the algebraic structure, the ordering
structure will provide a new powerful tool. This important idea has been widely used in solving integral
equations \cite{banas,carl,dhage1,amor,guo,aref,li,neto}.
In this work, we are mainly concerned with the existence results of solutions of the following system of
Quadratic nonlinear functional differential equation (in short QFDE)
\begin{equation}\label{fie.}
\displaystyle
\left\{
\begin{array}{llll}
\displaystyle\left(\frac{x(t)-f_1(t,x(t))}{f_2(t,y(t))}\right)'+ \lambda
\left(\frac{x(t)-f_1(t,x(t))}{f_2(t,y(t))}\right)
= g(t,y(t))\\
 y(t) = \displaystyle\frac{1}{1-b(t)|x(t)|}-p\left(t,\frac{1}{1-b(t)|x(t)|}\right)+p(t,y(t))\\ \\
\big(x(0),y(0)\big)=(x_0,y_0) \in \R^2,
\end{array}\right.
\end{equation}
where $\lambda \in \R_+$, $b:J\longrightarrow \R$ and $f_1,f_2,g,p:J\times \R\longrightarrow \R$ are continuous with $f_2$ not ever vanishing. The first equation of FDE \eqref{fie.} is general in the sense that it includes some important classes of functional
differential equations. If we take $x=y$, $f_1(t,x)=\phi(t,x)=0$ and $f_2(t,x) = f(t,x,y)$ and $g(t,x) =g(t,x,y)$, then  the FDE \eqref{fie.}
reduces to the following functional differential equation with a delay:
\begin{equation*}
\displaystyle
\left\{
\begin{array}{llll}
\displaystyle\left(\frac{x(t)}{f(t,x(t),X(t))}\right)'+ \lambda
\left(\frac{x(t)}{f(t,x(t),X(t))}\right)
= g(t,x(t),x_t),~~ t\in I
\\\\
x_0= \phi\in C([-r,0],\R), \ \ r>0
\end{array}\right.
\end{equation*}
whenever  $I=[0,T], \  T>0.$ The equation was examined in the paper \cite{dhage3,dhage4} and some special cases of this equation were
considered in \cite{mule}.\\
Again, when $f_1(t,x) =\phi(t,x)=0$, $f_2(t,x)=1$ and $g(t,x)=g(t,x,y)$, the QFDE \eqref{fie.} reduces to the following known nonlinear differential equation with maxima
\begin{equation*}
\displaystyle
\left\{
\begin{array}{llll}
 x'(t)+ \lambda x(t) = g(t,x(t),X(t),~~ t\in I\\\\
x(0)= \alpha_0 \in \R_+.
\end{array}\right.
\end{equation*}
The above nonlinear differential equation with maxima  has already been discussed in \cite{bainov} for
existence and uniqueness of the solutions via classical methods of Schauder and Banach fixed point principles.\\
Furthermore, if $f_1(t,x) = f(t,x,y)$, $f_2(t,x)=1$, $\phi(t,x)=0$ and $g(t,x)=g(t,x,y)$, the QFDE \eqref{fie.} reduces to the following HDE without maxima,
\begin{equation}
\displaystyle
\left\{
\begin{array}{llll}
(x(t)-f(t,x(t),X(t))'+ \lambda (x(t)-f_1(t,x(t)))= g(t,x(t),X(t))\\ \\
x(0))=\alpha_0 \in \R_+,
\end{array}\right.
\end{equation}
which is discussed in \cite{dhage0} via Dhage iteration method and established the existence and approximation result under some mixed partial Lipschitz and partial compactness conditions.\\

Note that the system \eqref{fie.} may be transformed into the following fixed point problem of the $2\times2$
block operator matrix
\begin{eqnarray}\label{d}
 \displaystyle
\left(
  \begin{array}{cc}
    A & B\cdot B' \\
    C & D \\
  \end{array}
\right),
\end{eqnarray}
where the entries of the matrix are, in general, nonlinear operators defined on partially ordered Banach
algebras. The operators occurring in the representation \eqref{d} are nonlinear, and our assumptions are as follows: $A$ a maps nondecreasing in a partially ordered Banach algebra $E$ into $E$ and $B,C,D$ and $B'$ are nondecreasing and positive operators from $E$ into $E$.\\

In this direction, the authors A. Jeribi, B. Krichen and N. Kaddachi in \cite{najib, najib1} have  established some fixed point for a $2\times 2$ operator matrix \eqref{d}, when $X$ is a Banach algebra satisfying certain condition. It is important to mention that the theoretical study was based on the existence of a solution of the following
equation
\begin{equation}\label{abc}
x=Ax\cdot Bx+Cx
\end{equation}
and obtained a lot of valuable results  (\cite{dhage6,dhage1, dhage2} and the references therein). These studies were
mainly based on the closure of the bounded domain, and properties of the operators $A$, $B$ and $C$ (cf.
partially completely continuous, partially nonlinear k-set contractive, partially condensing, and the potential
tool of the axiomatic partially measure of noncompactness,...).\\

This paper is organized as follows. In the next section, we give some preliminary results needed in the sequel.
In Section $3$, we present existence results for Equation \eqref{abc}. In Section $4$, we will deal with some
fixed point results for $2\times2$ block operator matrices in partially ordered Banach algebra. The main results
of this section are Theorems \ref{th1} and \ref{th2}. In Section $5$, we give an application showing the
existence of solutions of the system \eqref{fie.} in partially ordered Banach algebra.

\section{Auxiliary Results}
 Throughout this paper, let $(E,\preceq,\|\cdot\|)$ be a partially ordered Banach algebra with zero element $\theta$. Two elements $x$ and $y$ in  $E$ are called comparable if either the relation $x \preceq y$ or $y \preceq x$ holds. A non-empty subset $C$ of $E$ is called a chain or totally ordered set if all elements of $C$ are comparable. It is know that $E$ is regular if $\{x_n\}$ a nondecreasing (resp. nonincreasing) sequence in $E$
 and $x_n\to x^*$ as $n\to\infty$, then $x_n\preceq x^*$ (resp. $x_n\succeq x^*)$ for
 all $n\in \N$. The conditions guaranteeing the regularity of $E$ many be found in Guo and Lakshmikentham
 \cite{guolak} and Nieto and Lopez \cite{neto} and the references therein.\\

 At the beginning of this section, we present some basic facts concerning the partially measures of
noncompactness in $E$. If $C$ is a chain in $E$, then $C'$ denotes the set of all limit points of $C$ in $E$.
The symbol $\overline{C}$ stands for the closure of $C$ in $E$ defined by $\overline{C}= C\cup C'$. The set
$\overline{C}$ is called a closed chain in $E$. Thus, $\overline{C}$ is the intersection of all closed chains
containing $C$.
 In what follows, we denote by $\mathcal{P}_{cl}(E)$, $\mathcal{P}_{bd}(E)$, $\mathcal{P}_{rcp}(E)$,
 $\mathcal{P}_{ch}(E)$, $\mathcal{P}_{bd,ch}(E)$, $\mathcal{P}_{rcp,ch}(E)$, the class of all nonempty and
 closed, bounded, relatively compact, chains, bounded chains and relatively compact chains of $E$ respectively.
 Recall that the notion of the partial Kuratowski measure of noncompactness $\alpha^p(.)$ on $E$
 by the formula:
$$\alpha^p(C) =\displaystyle\inf\left\{r>0 , C=\bigcup_{i=1}^{n}C_i, \textrm{diam}(C_i)\leq r \ \forall i\right\}$$
where  diam$(C_i)=\sup\{\| x-y\| : x,y\in C_i\}$. For convenience we recall some basic properties of $\alpha^p(.)$ needed below \cite{dhage,dhage1,dhage2,dhage3}.
\begin{df}
A mapping $\alpha^p :\mathcal{P}_{bd,ch}(E)\longrightarrow \R_+$ is said to be a partially measure of noncompactness
in $E$ if it satisfies the following conditions:
\begin{enumerate}
  \item $\varnothing \neq (\alpha^p)^{-1}(\{0\})\subset \mathcal{P}_{rcp,ch}(E)$,
  \item $\alpha^p(\overline{C}) = \alpha^p(C)$,
  \item $\alpha^p$ is nondecreasing,
  \item If $\{C_n\}$ is a sequence of nondecreasing closed chains from $\mathcal{P}_{bd,ch}(E)$ with
      $\displaystyle\lim_{n\to \infty}\alpha^p(C_n)=0,$ then $\overline{C}_\infty=\displaystyle\cap_{n=0}^{\infty}C_n$ is a nonempty set and $\alpha^p(C_\infty)=0$,\\

The family of sets described in $1.$ is said to be the kernel of the measure of noncompactness $\alpha^p$ and is defined as
$$\ker\alpha^p=\{C\in \mathcal{P}_{bd,ch}(E)|\alpha^p(C)=0\}$$ \\
Clearly, $\ker\alpha^p\subset\mathcal{P}_{rcp,ch}(E)$. Observe that the intersection set $C_\infty$ from condition $4.$ is a member of the family $ker\alpha^p$. In fact, since $ \alpha^p (C_\infty) \leq \alpha^p (C_n)$ for any $n$, we infer that $\alpha^p (C_\infty) = 0.$ This yields that $ \alpha^p (C_\infty) \in \ker\alpha^p$. This simple observation will be essential in our further investigations.\\

The partially measure $\alpha^p$  of noncompactness is called sublinear if it satisfies

  \item $\alpha^p(C_1 + C_2) \leq \alpha^p(C_1) + \alpha^p(C_2),$ for all $C_1,C_2 \in \mathcal{P}_{bd,ch}$,
  \item $\alpha^p(\lambda C_1) = |\lambda|\alpha(C_1)$, for all $\lambda \in\R,$\\

  Again, $\alpha^p$ is said to satisfy maximum property if

  \item $\alpha^p(C_1\cup C_2)= \max\{\alpha^p(C_1),\alpha^p(C_2)\}$. \\

  Finally, $\alpha^p$ is said to be full or complete if
 \item  $\ker\alpha^p=\mathcal{P}_{rcp,ch}(E)$ $\hfill\diamondsuit$
\end{enumerate}
\end{df}
The following definitions (see \cite{dhage,dhage1,dhage2} and the references therein) are frequently used in the subsequent part of this paper.

\begin{df}
A mapping $T:E \longrightarrow E$ is called monotone nondecreasing if it preserves the order relation $\preceq$,
that is, if $x \preceq y$ implies $Tx \preceq Ty$ for all $x, y \in E$. Similarly, T is called monotone
nonincreasing if $x \preceq y $ implies $T x\succeq T y$ for all $x, y \in E$. A monotone mapping $T$ is one
which is either monotone nondecreasing or monotone nonincreasing on $E$.$\hfill\diamondsuit$
\end{df}

\begin{df}
A mapping $\psi:\R_+\longrightarrow\R_+$ is called a dominating function or, in short, $\mathcal{D}$-function if it is an upper semi-continuous and monotonic nondecreasing function satisfying $\psi(0)=0$$\hfill\diamondsuit$
\end{df}

\begin{df}
A mapping $T:E\longrightarrow E$ is called  partially nonlinear $\mathcal{D}$-Lipschitzian if there exist a $\mathcal{D}$-function $\psi:\R_+\longrightarrow \R_+$ satisfying
$$\|Tx-Ty\|\leq\psi(\|x-y\|),$$
for all comparable elements $x,y \in E$ where $\psi(0) = 0$. The function $\psi$ is called a
$\mathcal{D}$-function of $T$ on $E$. If $\psi(r)=kr$, $k>0$, then $T$ is called partially Lipschitzian with the Lipschitz constant $k$. In particular, if $k< 1$, then $T$ is called a partially contraction on $E$ with the contraction constant $k$. Finally, $T$ is called a partially nonlinear $\mathcal{D}$-contraction if it is a partially nonlinear $\mathcal{D}$-Lipschitzian with  $\psi(r)<r$ for $r>0$.$\hfill\diamondsuit$
\end{df}
\begin{rem}
Obviously, every partially Lipshitzian mapping is partially nonlinear $\mathcal{D}$-Lipshitizian. the converse may be not true.
\end{rem}
\begin{df}
A nondecreasing mapping $T : E \longrightarrow E$ is called partially nonlinear $\mathcal{D}$-set-Lipschitzian
if there exists a $\mathcal{D}$-function $\psi$ such that
 $$\alpha^p(TC)\leq \psi(\alpha^p(C)),$$
for all bounded chain $C$ in $E$. $T$ is called partially $k$-set-Lipschitzian if $\psi(r)=kr$, $k>0$. $T$ is called partially $k$-set-contraction if it is a partially $k$-set-Lipschitzian with $k < 1$. Finally, $T$ is called a partially nonlinear $\mathcal{D}$-set-contraction if it is a partially nonlinear
$\mathcal{D}$-Lipschitzian with $\psi(r)<r$ for $r>0$.$\hfill\diamondsuit$
\end{df}
\begin{df}
A mapping $T :E \longrightarrow E$ is called partially continuous at a point $a \in E$ if for $\varepsilon >0$ there exist a $\delta>0$ such that $\|Tx-Ta\|<\varepsilon$ whenever $x$ is comparable to $a$ and
$\|x-a\|<\delta$. $T$ called partially continuous on $E$ if it is partially continuous at every point of it. It is clear that if $T$ is partially continuous on $E$, then it is continuous on every chain $C$ contained in $E$.$\hfill\diamondsuit$
\end{df}
\begin{df}
A mapping $T : E\longrightarrow E $ is called partially bounded if $T(C)$ is bounded for every chain $C$ in $E$.
$T$ is called uniformly partially bounded if all chains $T (C)$ in $E$ are bounded by a unique constant. $T$ is
called bounded if $T(E)$ is a bounded subset of $E$.$\hfill\diamondsuit$
\end{df}
\begin{df}
A mapping $T : E \longrightarrow E$ is called partially compact if $T(C)$ is a relatively compact subset of $E$
for all totally ordered sets or chains $C$ in $E$. $T$ is called uniformly partially compact if $T(C)$ is a
uniformly partially bounded and partially compact on $E$. $T$ is called partially totally bounded if for any
totally ordered and bounded subset $C$ of $E$, $T(C)$ is a relatively compact subset of $E$. If $T$ is partially
continuous and partially totally bounded, then it is called partially completely continuous on $E$.
$\hfill\diamondsuit$
\end{df}
\begin{rem}
Note that every compact mapping on a partially normed linear space is partially compact and every partially
compactmapping is partially totally bounded, however the reverse implications do not hold. Again, every
completely continuous mapping is partially completely continuous and every partially completely continuous
mapping is partially continuous and partially totally bounded, but the converse may not be true. $\hfill\diamondsuit$
\end{rem}
\begin{df}
The order relation $\preceq$ and the metric $d$ on a non-empty set $E$ are said to be compatible if $\{x_n\}$ is
a monotone, that is, monotone nondecreasing or monotone nondecreasing sequence in E and if a subsequence
$\{x_nk\}$ of $\{x_n\}$ converges to $x^*$ implies that the whole sequence $\{x_n\}$ to $x^*$. Similarly, given
a partially ordered normed linear space $(E,\preceq,\|\cdot\|)$ the order relation $\preceq$ and the norm
$\|\cdot\|$ are said to be compatible if $\preceq$ and the metric $d$ defined through the norm $\|\cdot\|$ are
compatible.$\hfill\diamondsuit$
\end{df}
\begin{df}
A map $T: E\longrightarrow E$ is called $T$-orbitally continuous on $E$ if for any sequence $\{x_n\}\subseteq \mathcal{O}(x;T)=\{x,Tx,T^2x,\ldots,T^nx,\ldots\}$ we have that $x_n\to x^*$ implies that $Tx_n\to Tx^*$ for each $x\in E$. The metric space $E$ is called $T$-orbitally complete if every cauchy sequence $\{x_n\}\subseteq \mathcal{O}(x;T)$ converses to a point $x^*\in E$.
\end{df}
We need the following results in the sequel.
Let $(E,\preceq,\parallel\cdot\parallel)$ be partially ordered Banach algebra. Denote
$$E^+ =\{x\in E/x\succeq\theta\} \textrm{ and } \mathcal{K}=\{E^+\subset E/uv\in E^+ \textrm{ for all }u,v \in E^+\},$$
where $\theta$ is the zero element of $E$. The members $\mathcal{K}$ are called positive vectors in $E$.
\begin{lem} \textup{\cite{dhage3}}
If $u_1,u_2,v_1,v_2\in \mathcal{K}$ are such that $u_1\preceq v_1$ and $u_2\preceq v_2$, then $u_1u_2\preceq v_1v_2$.$\hfill\diamondsuit$
\end{lem}
\begin{df}
An operator $T:E\longrightarrow E$ is said to be positive if the range $R(T)$ of $T$is such that $R(T)\subseteq\mathcal{K}$.$\hfill\diamondsuit$
\end{df}
\begin{lem} \textup{\cite{dhage1}}\label{dh}
If $C_1$ and $C_2$ are two bounded chains in a partially ordered Banach algebra $E$, then
$$\alpha^p(C_1\cdot C_2)\leq \|C_2\|\alpha^p(C_1)+\|C_1\|\alpha^p(C_2)$$
where $\|C\| = \sup \{\|c\|, c\in C\}$. $\hfill\diamondsuit$
\end{lem}
\begin{thm}\textup{\cite{dhage}}\label{lem1}
Let $(E,\preceq , \|.\|)$ be a partially ordered set and let $T : E \rightarrow E$ be a
nondecreasing mapping. Suppose that there is a metric $d$ in X such that $(E,d)$ is a
$T$-orbitally complete metric space. Assume that there exists a $\mathcal{D}$-function $\psi$ such that
$$d(Tx,Ty)\leq \psi(d(x,y))$$
for all comparable elements $x,y \in E$ satisfying $\psi(r)<r $ for $r>0$.Further assume that
either $T$ is $T$-orbitally continuous on $E$ or $E$ is such that if $\left\{x_n\right\}$ is a nondecreasing
sequence with $x_n \rightarrow \overline{x} \in E $, then $ x \preceq \overline{x}$ for all $n \in \mathbb{N}$.
If there is an element $x_0 \in E$ satisfying $x_0 \preceq Tx_0$ or $x_0\succeq Tx_0$, then $T$ has a fixed point which is further
unique if "every pair of elements in $E$ has a lower and an upper bound".$\hfill\diamondsuit$
\end{thm}
\begin{thm}\textup{\cite{dhage}}\label{dh}
Let $(E,\preceq ,\|\cdot\|)$ be a partially ordered Banach algebra. Let $T :E\longrightarrow E$ be a nondecreasing, partially compact and continuous mapping. Further if the order relation $\preceq$ and the norm $\|\cdot\|$ in $E$ are compatible and if there is an element $x_0 \in E$ satisfying $x_0 \preceq Tx_0$ or $x_0\succeq Tx_0$, then $T$ has a fixed point.$\hfill\diamondsuit$
\end{thm}
\begin{thm}\textup{\cite{dhage1}}\label{dh1}
Let $(E,\preceq ,\|\cdot\|)$ be a regular partially ordered complete normed linear space such that the order relation and $\preceq $ the norm $\|\cdot\|$ are compatible. Let $T : E \longrightarrow E$ be a nondecreasing,
partially continuous and partially bounded mapping. If $T$ is partially nonlinear $\mathcal{D}$-set-contraction and if there exists an element $x_0 \in E$ such that $x_0 \preceq Tx_0$ or $x_0\succeq Tx_0$, then $T$ has a fixed point $x^*$ and
the sequence $\{T^nx_0\}$ of successive iterations converges monotonically to $x^*$.$\hfill\diamondsuit$
\end{thm}
\section{Reformulation of Dhage's fixed point theorems}
In this section, we prove some fixed point theorems in  partially ordered Banach algebras. Our results are formulated using some Dhage's type fixed point. Now, we establish our first fixed point theorem by modifying some assumptions of Theorem 4.1 in \cite{dhage1}.
\begin{thm}\label{t1}
Let $S$ be non-empty closed and partially bounded subset of a regular partially ordered Banach algebra $(E,\preceq ,\|\cdot\|)$  such that the order relation $\preceq$ and the norm $\|\cdot\|$ in E are compatible. Let $A:E\longrightarrow\mathcal{K}$, $B:S\longrightarrow\mathcal{K}$ and $C:E\longrightarrow E$ be three nondecreasing operators satisfying the following conditions:\\
\noindent $(i)$ $A$ and $C$ are  partially nonlinear $\mathcal{D}$-Lipschitzians with $\mathcal{D}$-functions $\psi_A$ and $\psi_ C,$\\
\noindent $(ii)$ $B$ is  continuous and partially compact,\\
\noindent $(iii)$ There exists  $x_0\in S$ such that $x_0\preceq Ax_0\cdot By + Cx_0$ or $x_0\succeq Ax_0\cdot By + Cx_0$ for each $y\in S$,\\
\noindent $(iv)$ $Ax\cdot By + Cx \in S$ for all $y\in S$\\
\noindent $(v)$ Every pair of elements $x,y \in E$ has a lower and an upper bound in $E$.\\
\noindent Then the equation \eqref{abc} has a fixed point in $E$as soon as $M\psi_A(r) + \psi_C(r) < r$ if $r>0$, where $M=\|B(E)\|$. $\hfill\diamondsuit$
\end{thm}
\begin{pf}
Suppose that there exists $x_0\in S$  such that $x_0\preceq Ax_0\cdot By+Cx_0$.
Let $y \in S$ and define a mapping define a mapping
$A_y:S\longrightarrow S$ by the formula
$$A_y(x)=Ax\cdot By+Cx.$$
Since $A$ and $B$ are positive and $A$, $B$ and $C$ are nondecreasing, $A_y$ is nondecreasing
Now, let $x_1, x_2 \in E$ be two comparable elements. Then
\begin{eqnarray*}
\|A_y(x_1)-A_y(x_2)\| &=&\|Ax_1\cdot By+Cx_1 -Ax_2\cdot By-Cx_2\|\\
                       &\leq& \|Ax_1\cdot By -Ax_2\cdot By\|+\|Cx_1 -Cx_2\|\\
                       &\leq& \|By\|\|Ax_1 -Ax_2\|+\|Cx_1 -Cx_2\|\\
                       &\leq& \left(M\psi_A+\psi_C\right)(\|x_1-x_2\|).
\end{eqnarray*}
This implies that the operator $A_y$ is a partially nonlinear $\mathcal{D}$-contraction. Hence,  by Theorem \ref{lem1} there exist a unique point $x^*\in E$ such that
$$x^*= Ax^*\cdot By+Cx^*.$$
because the hypothesis $(iv)$ for all $y\in S$ we have $x^*\in S$. Define a mapping
\begin{eqnarray*}
T:&S&\longrightarrow S\\
  &y&\longmapsto x^*,
\end{eqnarray*}
where $x^*$ is the unique solution of the equation $Ax^*\cdot By+Cx^*$. Since $A, B$ and $B'$ are nondecreasing and $B$ and $B'$ are positive , $T$ is nondecreasing.  Now we show that $T$ is continuous. Let  $\{y_n\}_{n=0}^{\infty}$ be any sequence in $B(E)$ converging to a point $y$, Since $S$ is closed, $y\in S$.
Then,
\begin{eqnarray*}
\|Ty_n - Ty\| &=& \|Ax^*_n\cdot y_n + Cx^*_n - Ax\cdot y - Cx\|\\
            &\leq&  \|Ax^*_n\cdot y_n - Ax\cdot y\| + \|Cx^*_n- Cx\|\\
            &\leq&  \|Ax^*_n\cdot y_n - Ax\cdot y_n\| + \|Ax\cdot y_n - Ax\cdot y\| + \|Cx^*_n- Cx\| \\
            &\leq&  (M\psi_A+\psi_c)(\|x_n- x\|)+ \|Ax\|\|y_n -y\|.
\end{eqnarray*}
Hence
$$\lim\sup_{n}\|Ty_n-Ty\|\leq  (M\psi_A+\psi_c)(\lim\sup_{n}\|x_n- x\|)+ \|Ax\|\lim\sup_{n}\|y_n -y\|.$$
This show that $\lim\limits_{n \rightarrow +\infty} \|Ty_n - Ty\|=0$ and  the claim is
approved.
Next, we shows that $T$ is partially compact. In fact, let $\mathcal{C}$ be a chain in $S$,
for any $z\in \mathcal{C}$ we have
\begin{eqnarray*}
\|Az\|&\leq&\|Aa\|+\psi_A(\|z-a\|)\\
      &\leq&\|A_a\|+\frac{\|z-a\|}{M}\\
      &\leq& c,
\end{eqnarray*}
where $c=\|A_a\|+\displaystyle\frac{diam \mathcal{C}}{M}$ for some fixed point $a$ in $\mathcal{C}$.\\
Let $\varepsilon>0$ be given. By assumption $(ii)$, we infer that $B(\mathcal{C})$ is partially
totally bounded, then there exist a chain $Y=\{y_1,\ldots,y_n\}$ of point in $\mathcal{C}$
such that
$$B(\mathcal{C}) \subseteq \bigcup_{i=1}^{n}\mathbf{B}_\delta (w_i),$$
where $w_i=By_i$ and $\delta= \frac{1}{c}(\varepsilon-(M\phi_A(\varepsilon)+\phi_C(\varepsilon)))$ and $\mathbf{B}_\delta (w_i)$ is an open ball in $E$ centered at
$w_i$ of radius $\delta$. Therefore, for any $y\in\mathcal{C}$, we have $y_k\in Y$ such that
$$c\|By_i-By\|\leq \delta.$$
Also, we have
\begin{eqnarray*}
\|Ty_k-Ty\| &\leq&\|Ax^*_k\cdot By_k-Ax\cdot By\|-\|Cx^*_k-Cx\|\\
                &\leq& M(\psi_A+\psi_C)\|x^*_k-x\|+\|Ax\|\|By_k-By\|\\
                &\leq& M(\psi_A+\psi_C)\|x^*_k-x\|+c\|By_k-By\|.\\
\end{eqnarray*}
Then
\begin{eqnarray*}
(I-M(\psi_A+\psi_C))(\|Ty_k-Ty\|) \leq c\|By_k-By\|.
\end{eqnarray*}
So,
$$\|T(y_k)-T(y)\|<\varepsilon,$$
because $y\in \mathcal{C}$ is arbitrary,
$$T(S)\subseteq \displaystyle\bigcup_{i=1}^{n} \mathbf{B}_\varepsilon(k_i),$$
where $k_i = T(y_i)$. As a result, $T(S)$ is partially totally bounded in $E$. Hence, $T$ is is partially compact. The order relation $\preceq$ and the norm $\|.\|$ are compatible, so the desired conclusion follows of Theorem  $\ref{dh}$ we have $T$ has a fixed point in $S$.
\end{pf}\\
Now, modifying same assumptions of  Theorem 4.8 in \cite{dhage2}, we have the following result.
\begin{thm}\label{t2}
Let $S$ be closed and partially bounded subset of a regular partially ordered Banach algebra $(E,\preceq ,\|\cdot\|)$ such that the order relation $\preceq$ and the norm $\|\cdot\|$ are compatible in every compact chain $C$ of $S$. Let $A,B:S\rightarrow \mathcal{K}$ and $C:S\longrightarrow E$ be three nondecreasing  operators satisfying the following conditions\\
\noindent $(i)$ $A$ and $C$ are partially nonlinear $\mathcal{D}$-Lipschitzians with $\mathcal{D}$-functions $\psi_A$ and $\psi_ C$ respectively,\\
\noindent $(ii)$ $B$ is partially completely continuous ,\\
\noindent $(iii)$ $\displaystyle\left(\frac{I-C}{A}\right)^{-1}$ exist on $B(S)$ and is nondecreasing,\\
\noindent $(iv)$ There exist $x_0\in S$ such that $x_0\preceq Ax_0\cdot By + Cx_0$ or $x_0\succeq Ax_0\cdot B y + Cx_0,$ for each $y\in S$, \\
\noindent $(v)$ $Ax\cdot By+ Cx\in S$ for all $y \in S,$\\
\noindent Then the equation \eqref{abc} has a fixed point in $S$ as soon as $M\psi_A(r) + \psi_C(r) < r$ if $r>0$ where
$M=\|B(S)\|$.
\end{thm}
\begin{pf}
Suppose that there exist $x_0\preceq Ax_0\cdot By+Cx_0$.  It is easy to check  that the vector $x\in E$ is a solution for the equation \eqref{abc}, if and only if $x$ is a fixed point for the operator $T:=\displaystyle\left(\frac{I-C}{A}\right)^{-1}B.$ From assumption $(iii)$ it follows that, for each $y\in S$ there is a unique $x_y\in E$ such that
$$\left(\frac{I-C}{A}\right)x_y=By$$
or, in an equivalently way
$$Ax_y\cdot By+Cx_y=x_y.$$
Since the assumption $(v)$ holds, then $x_y\in S$. Hence the map $T:S\longrightarrow S$ is well define. from assumption $(iii)$, it follows that $T$ is nondecreasing.  Now, in view of Theorem \ref{dh}, it suffices to prove that $T$ is continuous and partially compact. Indeed
 let $\{y_n\}_{n=0}^{\infty}$ be any sequence in $S$ converging to a point $y$. Because $S$ is closed,
$y\in S$. Now
\begin{eqnarray*}
\|Ty_n - Ty\| &=&     \|Ax^*_n\cdot By_n + Cx^*_n - Ax^*\cdot By - Cx^*\|\\
            &\leq&  \|Ax^*_n\cdot By_n - Ax^*\cdot By\| + \|Cx^*_n- Cx^*\|\\
            &\leq&  \|Ax^*_n\cdot By_n - Ax^*\cdot By_n\| + \|Ax^*\cdot By_n - Ax^*\cdot By\| + \|Cx^*_n- Cx^*\|
            \\
            &\leq&  (M\psi_A+\psi_C)(\|x^*_n- x^*\|)+ \|ATy\|\|By_n -By\|.
\end{eqnarray*}
Hence,
$$\displaystyle\lim_n\sup \|Ty_n - Ty\|\leq (M\psi_A+\psi_C)\lim_n\sup\|Ty_n-Ty\|+\|ATy\|\lim_n\sup\|By_n-By\|.$$
This shows that $\displaystyle\lim_{n\to\infty}\|Ny_n - Ny\| = 0$ and the claim is approved.
Next, we shows that $T$ is paritally compact. Since $B$ is partially compact and  $\displaystyle\left(\frac{I-C}{A}\right)^{-1}$ is continuous, then the composition mapping  $T= \displaystyle\left(\frac{I-C}{A}\right)^{-1}B$  is continuous and partially compact on $E$.
Next, the order relation $\preceq$ and the norm $\|\cdot\|$ in $E$ are compatible.  Hence, an application of Theorem \ref{dh} infer that $T$ has, at least, one fixed point in $S$.
\end{pf}
%
%

\section{Application of Dhage's fixed point to block matrix operator}
In what follows, we will study the existence of a fixed point for the block matrix operators.
\begin{thm}\label{th1}
Let $(E,\preceq ,\|\cdot\|)$ be a regular partially ordered Banach algebra such that the order relation $\preceq$ and the norm $\|\cdot\|$ in E are compatible. Let $A,C,D :E\longrightarrow E$ and $B,B':E\longrightarrow \mathcal{K}$ be five  nondecreasing  operators satisfying the following assumptions:\\
\noindent $(i)$ $A$, $B$ and $C$ are  partially bounded and  partially nonlinear $\mathcal{D}$-Lipschitzians with $\mathcal{D}$-functions $\psi_A$, $\psi_B$ and $\psi_ C$ respectively,\\
\noindent$(ii)$ $(I-D)^{-1}$ exist and  partially nonlinear $\mathcal{D}$-Lipschitz
with $\mathcal{D}$-function $\psi_\phi$ and $(I-D)^{-1}C$ is nondecreasing,\\
\noindent $(iii)$ $B'$ is partially continuous and $C$ is compact,\\
\noindent $(iv)$ There exist $x_0\in E,$ such that $x_0 \preceq Ax_0 + Tx_0\cdot T'x_0$ or
$x_0\succeq Ax_0 + Tx_0\cdot T'x_0,$ where $T=B(I-D)^{-1}C$ and $ T'=B'(I-D)^{-1}C$.\\
Then the operator matrix \eqref{d} has a fixed point in $E\times E$, whenever $M\psi(r)\leq r $ if $r>0$ with $\psi(r)=\psi_B\circ \psi_\phi\circ\psi_C(r) + \psi_A(r)$ and $M=\|T'(E)\|.$
\end{thm}
\begin{pf}
Suppose that there is an element $x_0 \in E$ such that $x_0\preceq Ax_0+Tx_0\cdot T'x_0$.
Define a mapping $F:E\longrightarrow E$ by the formula
         $$Fx= Ax +Tx\cdot T'x.$$
Since $B$ and $B'$ are positive and $A, B, B'$ and $(I-D)^{-1}C$ are nondecreasing we infer that $F$ is nondecreasing.
 Next, we claim that $F$ is partially continuous. To do us, let $\{x_n; n\in \N\} $ be a sequence in $E$ which converge to a point $x$ such that $x_n$ and $x$ are comparable. From assumption $(i)$, it follows that:
\begin{eqnarray*}
  \|Fx_n-Fx\| &=& \|Ax_n + Tx_n\cdot T'x_n - Ax -Tx\cdot T'x\| \\
            &\leq&\|Ax_n-Ax\| + \|Tx_n\cdot T'x_n - Tx\cdot T'x\| \\
            &\leq&\|Ax_n -Ax\| +\|Tx_n\|\|T'x_n - T'_x\|+ \|T'x\|\|Tx_n-Tx\|\\
            &\leq& \psi_A +M\psi_B\circ \psi_\phi\circ \psi_C(\|x_n-x\|)+ \|Tx_n\|\|T'x_n - T'x\|.
\end{eqnarray*}
From the partially continuity of $B'$, and taking limit supremum in the aforementioned inequality yields that
$$\lim_{n\to \infty}\|Fx_n-Fx\|=0.$$
This proves that $F$ is a partially continuous operators on $E$.
 Again by assumption $(iv)$, it is clear that $x_0\preceq Fx_0$. Moreover, it easy to show that $F$ is partially bounded. Furthermore, we show that $F$ is a partially nonlinear $\mathcal{D}$-set-contraction on $E$. Let $\mathcal{C}$ be a bounded chain in $E$. Then by definition of $F$, we have
$$ F(\mathcal{C})\subseteq A(\mathcal{C})+ T(\mathcal{C})T'(\mathcal{C}).$$
Since $F$ is nondecreasing and partially continuous $F(\mathcal{C})$ is again a bounded chain in E. Keeping in mind the relatively compactness of $T'(\mathcal{C})$ and making use of Lemma \ref{dh} together with the subadditivity of the partially Kuratowski measure of noncompactness $\alpha^{p}$, enables us to have
\begin{eqnarray*}
  \alpha^{p}(F(\mathcal{C})) &\leq& \|T(\mathcal{C}))\|\alpha^{p}(T'(\mathcal{C}))+
  \|T'(\mathcal{C})\|\alpha^{p}(T(\mathcal{C}))+\alpha^{p}(A(\mathcal{C})) \\
                            &\leq& M\alpha^{p}(T(\mathcal{C})) + \alpha^{p}(A(\mathcal{C})) \\
                            &\leq& M\psi_B\circ \psi_\phi\circ \psi_C(\alpha^{p}(\mathcal{C}))+
                            \psi_A(\alpha^{p}(\mathcal{C})) \\
                            &=& \psi(\alpha^{p}(\mathcal{C})),
\end{eqnarray*}
where $\psi(r) = \psi_A(r)+M\psi_B\circ \psi_\phi\circ \psi_C(r)  < r $ for all $r>0.$
 This shows that $F$ is a partially nonlinear $\mathcal{D}$-set-contraction on $E$. Finally  the relation
 $\preceq$ and the norm $\|\cdot\|$ are compatible, so the desired conclusion follows by an application of Theorem  \ref{dh1} we have $x= Ax + Tx\cdot T'x$ has a solution in $E$. Now, use the vector $y=(I-D)^{-1}Cx$ to achieve the proof.
\end{pf}
\\
\begin{rem}\label{rq}
If $A$ is a contraction on $S$ into itself with constant contraction $k$ and $B$ a partially nonlinear $\mathcal{D}$-lipshitzian with a $\mathcal{D}$-function $\psi_B$, and $C(S)\subset (I-D)(E)$, then the operator inverse  $(I-A)^{-1}B$  exist and  partially nonlinear $\mathcal{D}$-lipshiztian with a $\mathcal{D}$-function $\psi(r) = \frac{1}{1-k}\psi_B(r)$.
\end{rem}
Next, we will combine Theorem \ref{th1} and Remark \ref{rq} in order to to obtain the following
fixed point theorem:
\begin{cor}
Let $(E,\preceq ,\|\cdot\|)$ be a regular partially ordered Banach algebra such that the order relation $\preceq$ and the norm $\|\cdot\|$ in $E$ are compatible. Let $A,C,D :E\longrightarrow E$ and $B,B':E\longrightarrow \mathcal{K}$ are five  nondecreasing  operators satisfying the following assumptions:\\
 \noindent $(i)$ $A$, $B$ and $C$ are partially bounded and  partially nonlinear $\mathcal{D}$-Lipschitzians with $\mathcal{D}$-functions $\psi_A$, $\psi_B$ and $\psi_ C$ respectively,\\
 \noindent $(ii)$  $D$ is a contraction with a contraction constant k and $(I-D)^{-1}C$ is nondecreasing,\\
 \noindent$(iii)$ $B'$ is partially continuous and $C$ is partially  compact, and $C(E)\subset (I-D)(E)$\\
 \noindent$(iv)$ $x_0 \preceq Ax_0 + Tx_0\cdot T'x_0$ or
 $x_0\succeq Ax_0 + Tx_0\cdot T'x_0,$ for some $x_0\in E,$ where $T=B(I-D)^{-1}C$ and $ T'=B'(I-D)^{-1}C$.\\
 Then the operator matrix \eqref{d} has a fixed point in $E\times E$, whenever $M\psi(r)\leq r $ if $r>0$ with $\psi(r)=\psi_B\circ(\frac{1}{1-k}) \psi_C(r)$ + $\psi_A(r)$ and $M=\|T'(E)\|.$
\end{cor}
Based on the conditions $C\equiv 0_E$ and $D=1_E$ in which $0_E$ and $1_E$ represents respectively the zero and the unit element of the partially ordered Banach algebra $E$, we infer from Theorem \ref{th1} the following result:
\begin{cor}\cite{dhage1}
Let $(E,\preceq ,\|\cdot\|)$ be a regular partially ordered Banach algebra such that the order relation $\preceq$ and the norm $\|\cdot\|$ in E are compatible. Let $A:E\longrightarrow E$ and $B,B':E\longrightarrow \mathcal{K}$ be five  nondecreasing  operators satisfying the following assumptions:\\
\noindent $(i)$ $A$, $B$ and  are  partially bounded and  partially nonlinear $\mathcal{D}$-Lipschitzians with $\mathcal{D}$-functions $\psi_A$, $\psi_B$ respectively,\\
\noindent $(ii)$ $B'$ is partially continuous and compact,\\
\noindent $(ii)$ There exist $x_0\in E,$ such that $x_0 \preceq Ax_0 + Bx_0\cdot B'x_0$ or
$x_0\succeq Ax_0 + Bx_0\cdot B'x_0$.\\
Then the operator matrix \eqref{d} has a fixed point in $E\times E$, whenever $M\psi(r)\leq r $ if $r>0$ with $\psi(r)=\psi_B + \psi_A(r)$ and $M=\|B'(E)\|.$
\end{cor}
Next, we can modify some assumptions of Theorem \ref{t2} in order to study the  problem in block operator matrix.
\begin{thm}\label{th2}
Let $S$ be a non-empty closed partially bounded subset of a regular partially ordered Banach algebra $(E,\preceq ,\|\cdot\|)$ such that the order relation $\preceq$ and the norm $\|\cdot\|$ are compatible in every compact chain $\mathcal{C}$ of $S$. Let $A,C:S\longrightarrow E$ and $B,B' ,D:E\longrightarrow \mathcal{K}$ be five  nondecreasing  operators satisfying the following assumptions:\\
\noindent $(i)$ $A$, $B$ and $C$ are partially nonlinear $\mathcal{D}$-Lipschitzians
with $\mathcal{D}$-functions $\psi_A$, $\psi_B$ and $\psi_ C$ respectively,\\
\noindent$(ii)$  $(I-D)^{-1}$ is nondecreasing and  partially nonlinear $\mathcal{D}$-Lipschitzian with $\mathcal{D}$-functions $\psi_\phi$,\\
\noindent$(iii)$ $B'$ is continuous and $C$ is a partially compact operator such that $C(S)\subseteq (I-D)(S)$\\
\noindent $(iv)$ $Ax+Tx\cdot T'y \in S$ for all $y\in S$ where $T=B(I-D)^{-1}C$ and $ T'=B'(I-D)^{-1}C$,\\
\noindent$(v)$ $\displaystyle\left(\frac{I-A}{T}\right)^{-1}$ exist on $T'(S)$ and nondecreasing, \\
\noindent$(vi)$ $x_0 \preceq Ax_0 + Tx_0\cdot T'y$ or
$x_0\succeq Ax_0 + Tx_0\cdot T'y,$ for each $y\in S,$ \\
\noindent$(vii)$ every pair of elements $x,y \in E$ has a lower and an upper bound in $E$. \\
Then the operator matrix \eqref{d} has a fixed point in $S\times S$, whenever $M\psi(r)\leq r $ if $r>0$ with $\psi(r)=\psi_B\circ \psi_\phi\circ\psi_C(r) + \psi_A(r)$ and $M=\|T'(S)\|.$
\end{thm}
\begin{pf}
Suppose that there is an element $x_0 \in E$ such that $x_0\preceq Ax_0+Tx_0\cdot T'y$.
From assumption $(iv)$, it follows that for each $y\in S$, there exist a unique point $x\in S$ such that
$$\displaystyle\left(\frac{I-A}{T}\right)x=T'y$$
or, equivalently
$x=Ax+Tx\cdot T'y$. Because assumption $(iv)$ hold, then $x\in S$. Therefore, we can define
  $F:S\longrightarrow S$ by the formula
\begin{equation}
F(x) = \left(\frac{I-A}{T}\right)^{-1}T'x.
\end{equation}
In view of assumption $(v)$, it follows that, $F$ is nondecreasing on $E$. Next, we will prove that $F$ is continuous. To see that, let $\{z_n\}_{n=0}^{\infty}$ be any sequence on $E$ such that $z_n\to z$, and let
\begin{equation*}\label{b}
\displaystyle
\left\{
\begin{array}{llll}
y_n =T'(z_n) \textrm{  and  } y = T'(z)\\
x_n = \displaystyle\left(\frac{I-A}{T}\right)^{-1}(y_n)  \textrm{  and  }
x = \displaystyle\left(\frac{I-A}{T}\right)^{-1}(y).
\end{array}\right.
\end{equation*}
Then, it easy to show that $y_n\to y$, and we have
\begin{equation*}\label{b}
\displaystyle
\left\{
\begin{array}{llll}
x_n &=& Ax_n +Tx_n\cdot T'y_n \\
x &=& Ax +Tx\cdot T'y.
\end{array}\right.
\end{equation*}
So,
\begin{eqnarray*}
\|x_n - x\| &=&     \|Ax_n + Tx_n\cdot T'y_n - Ax + Tx\cdot T'y\|\\
            &\leq&  \|Ax_n - Ax\| + \|Tx_n\cdot T'y_n - Tx\cdot T'y\|\\
            &\leq&  \|Ax_n.y_n - Ax.y_n\| + \|Ax.y_n - Ax.y\| + \|Cx_n- Cx\| \\
            &\leq&  \psi_A + M\psi_B\circ\psi_\phi\circ\psi_C(\|x_n- x\|)+ \|Tx\|\|y_n -y\|,
\end{eqnarray*}
where $\psi(r) = \psi_A (r)+ (M\psi_B\circ\psi_\phi\circ\psi_C)(r)<r$ for all $r>0$ and the constant $M$ exists in view of the fact $T'$ is partially bounded operator on $S$. Taking limit supermum in the aforementioned inequality yields that $F$ is continuous. Further, from assumption $(iii)$, we show that $B'(I-D)^{-1}C$ as well as $F$ is partially compact. Finally the relation $\preceq$ and the norm $\|\cdot\|$ are compatible. Hence, an application of Theorem  \ref{dh} infer that $F$ has at least, one solution in $S\times S$. Now, use the vector $y=(I-D)^{-1}Cx$ to achieve the proof.
\end{pf}
\\
As a consequence we have the following fixed point result.
\begin{cor}\label{cor}
Let $S$ be a non-empty closed partially bounded subset of a regular partially ordered Banach algebra $(E,\preceq ,\|\cdot\|)$ such that the order relation $\preceq$ and the norm $\|\cdot\|$ in E are compatible. Let $A,C:S\longrightarrow E$ and $B,B' ,D:E\longrightarrow \mathcal{K}$ be five  nondecreasing  operators satisfying the following assumptions:\\
\noindent $(i)$ $A$, $B$ and $C$ are partially bounded and partially nonlinear $\mathcal{D}$-Lipschitzians
with $\mathcal{D}$-functions $\psi_A$, $\psi_B$ and $\psi_ C$ respectively,\\
\noindent$(ii)$  $D$ is contraction with a contraction constant $k$ and $C(S)\subset(I-D)(S)$,\\
\noindent$(iii)$ $B'$ is  continuous and $C$ is a partially compact operator,\\
\noindent $(iv)$ $x=Ax+Tx\cdot T'y$, $y\in S \Rightarrow x\in S$  where $T=B(I-D)^{-1}C$ and $ T'=B'(I-D)^{-1}C$,\\
\noindent$(v)$ $\displaystyle\left(\frac{I-A}{T}\right)^{-1}$ exists on $T'(S)$ and nondecreasing, \\
\noindent$(vi)$ $x_0 \preceq Ax_0 + Tx_0\cdot T'x_0$ or
$x_0\succeq Ax_0 + Tx_0\cdot T'x_0,$ for some $x_0\in E,$ \\
\noindent$(vii)$ every pair of elements $x,y \in E$ has a lower and an upper bound in $E$. \\
Then the operator matrix \eqref{d} has a fixed point in $S\times S$, whenever $M\psi(r)\leq r $ if $r>0$ with $\psi(r)=\psi_B\circ(\frac{1}{1-k}) \psi_C(r)$ + $\psi_A(r)$ and $M=\|T'(S)\|.$
\end{cor}
From Theorem \ref{th2} and corollary \ref{cor} , without any hurdle we can derive
the following corollary.
\begin{cor}
Let $S$ be a non-empty closed partially bounded subset of a regular partially ordered Banach algebra $(E,\preceq ,\|\cdot\|)$ such that the order relation $\preceq$ and the norm $\|\cdot\|$ in E are compatible. Let $A,C:S\longrightarrow E$ and $B,B' ,D:E\longrightarrow \mathcal{K}$ be five  nondecreasing  operators satisfying the following assumptions:\\
\noindent $(i)$ $A$, $B$ and $C$ are partially bounded and partially nonlinear $\mathcal{D}$-Lipschitzians
with $\mathcal{D}$-functions $\psi_A$, $\psi_B$ and $\psi_ C$ respectively,\\
\noindent$(ii)$  $D$ is  contraction with a contraction constant $k$ and $C(S)\subset(I-D)(E)$,\\
\noindent$(iii)$ $B'$ is partially completely continuous \\
\noindent $(iv)$ $x=Ax+Tx\cdot T'y$, $y\in S \Rightarrow x\in S$ where $T=B(I-D)^{-1}C$ and $ T'=B'(I-D)^{-1}C$,\\
\noindent$(v)$ $\displaystyle\left(\frac{I-A}{T}\right)^{-1}$ exist on $B'(E)$ and nondecreasing,\\
\noindent$(vi)$ $x_0 \preceq Ax_0 + Tx_0\cdot T'x_0$ or
$x_0\succeq Ax_0 + Tx_0\cdot T'x_0,$ for some $x_0\in E,$ \\
\noindent$(vii)$ every pair of elements $x,y \in E$ has a lower and an upper bound in $E$. \\
Then the operator matrix \eqref{d} has a fixed point in $S\times E$, whenever $M\psi(r)\leq r $ if $r>0$ with $\psi(r)=\psi_B\circ(\frac{1}{1-k}) \psi_C(r)$ + $\psi_A(r)$ and $M=\|T'(E)\|.$
\end{cor}
\section{Existence solutions for a  system of functional differential equation}
The FDE \eqref{fie.} is considered in the function space $E=C(J, \R)$ of continuous real-valued functions defined on $J$. We define  a norm $\|\cdot\|$ and the order relation $\preceq$ in $C(J,\R)$ by
\begin{equation}\label{1}
\displaystyle\|x\|_\infty=\sup_{t\in J}|x(t)|
\end{equation}
and
\begin{equation}\label{2}
x\leq y\Longleftrightarrow x(t)\leq y(t)
\end{equation}
for all $t\in J$. Clearly, $C(J,\R)$ is a Banach space with respect to above supremum norm and also partially ordered with respect to the above partially order relation $\leq$. It is known that the partially ordered Banach space $C(J,\R)$ has some nice properties with respect to the above order relation in it. The following lemma follows by an application of Arzela-Ascolli theorem.
\begin{lem}\cite{dhage sb}
Let $(C(J,\R),\leq,\|\cdot\|)$ be a partially ordered Banach space with the norm $\|\cdot\|$ and the order relation $\leq$ defined by \eqref{1} and \eqref{2} respectively. Then $\|\cdot\|$ and $\leq$ are compatible in every partially compact subset of $C(J,\R)$.
\end{lem}
\noindent The purpose of this section is to apply theorem \ref{th1} to discuss the existence of solutions for the following nonlinear quadratic functional differential  equations QFDE \eqref{fie.}.\\
\noindent We need the next definition in what follows.
\begin{df}
A functions $u,v \in C(J,\R)$ is said to be a lower solution of the \eqref{fie.}, if it satisfies
\begin{equation*}
\displaystyle
\left\{
\begin{array}{llll}
\displaystyle\left(\frac{u(t)-f_1(t,u(t))}{f_2(t,v(t))}\right)'+ \lambda
\left(\frac{u(t)-f_1(t,u(t))}{f_2(t,v(t))}\right)
\leq g(t,v(t))\\
\displaystyle v(t) \leq \displaystyle\frac{1}{1-b(t)|u(t)|}-p\left(t,\frac{1}{1-b(t)|u(t)|}\right)+p(t,v(t))\\ \\
\big(u(0),v(0)\big)\leq(u_0,v_0) \in \R^2.
\end{array}\right.
\end{equation*}
for all $t\in J$. Similarly, a functions $u',v' \in C(J,\R)$ is said to be an upper solution of the FDE \eqref{fie.} if it satisfies the conditions above with the inequality reversed.
\end{df}
\noindent We consider the following set of assumptions in what follows:\\

\noindent $(\mathcal{H}_0)$ The functions $b:J\longrightarrow\R$ is continuous.\\
\noindent $(\mathcal{H}_1)$ The function $t\longmapsto f_1(t,0)$ is bounded on $J$ with bound $F_0$.\\
\noindent $(\mathcal{H}_2)$ $f_2$ define a function   $f_2:J\times\R\longrightarrow\R_+$ are nondecreasing in $x$ for all $t\in J$.\\
\noindent $(\mathcal{H}_3)$  There exist two constants $L,K \in \R_+^*$ such that
$$0<f_i(t, x(t)) - f_i(t, y(t)) \leq \frac{L(x-y)}{K + (x-y)}$$
\hspace*{20pt}  for all $t\in J$ , $i=1,2$ and $x , y \in \R$ with $x\geq y$. Moreover, $L \leq K$.\\

\noindent $(\mathcal{H}_4)$  The function  $p:J\times\R\longrightarrow\R $ is nondecreasing in $x$ for all $t\in J$ and contraction with a constant $k$. \\


\noindent $(\mathcal{H}_5)$ The function $g:J\times\R\longrightarrow\R_+$ is nondecreasing and there exist a function $h\in L^1(J)$ such that
$$\|g(t,x)\|\leq \|h\|_{L^1} \textrm{ for all } t\in J \textrm{ and } x\in\R.$$

\noindent $(\mathcal{H}_6)$ The QFDE \eqref{fie.} has a lower solution in $C(J,\R)\times C(J,\R)$\\

 \begin{rem}
 Assume that the assumption $(\mathcal{H}_6)$ holds. Then a functions $x\in C(J,\R)$ is a solution of
 \eqref{fie.} if and only if is a solution of the functional  integral equation
 \begin{equation}\label{fe}
\displaystyle
\left\{
\begin{array}{llll}
 x(t)&=& f_1(t,x(t))+f_2(t,y(t))\left(ce^{-\lambda t}+e^{-\lambda t}\displaystyle\int_{0}^{t}e^{\lambda
s}g(s,y(s))ds\right)\\
\displaystyle y(t) &=& \displaystyle\frac{1}{1-b(t)|x(t)|}-p\left(t,\frac{1}{1-b(t)|x(t)|}\right)+p(t,y(t))
\end{array}
\right.
\end{equation}
where $\displaystyle c = \frac{x_0- f_1(0,x_0)}{f_2(0,y_0)}.$
\end{rem}
Indeed, let $g\in C(J,\R)$. Assume first that $x$ is a solution of the QFDE \eqref{fie.} define on $J$. By definition, the function $t\longmapsto\displaystyle\frac{x(t)-f_1(t,x(t))}{f_2(t,y(t))}$ is continuous on $J$, and so, differentiable there, whence $\left(\displaystyle\frac{x(t)-f_1(t,x(t))}{f_2(t,y(t))}\right)'$ is integrable on $J$. Applying integration to \eqref{fie.} from $0$ to $t$, we obtain the HIE \eqref{fe} on $J$.\\

Conversely, assume that $x$ satisfies Eq. \eqref{fe}. Then by direct differentiation we obtain the first
equation in QFDE \eqref{fie.}. Again, substituting $t = 0$ in Eq. \eqref{fe} yields
$$\displaystyle\frac{x(0)-f_1(0,x(0))}{f_2(0,y(0))}=\frac{x_0-f_1(0,x_0)}{f_2(0,y_0)}$$
whence, $(x(0),y(0))=(x_0,y_0)$. Thus, the QFDE \eqref{fie.} holds.\\


\noindent Now, we are in a position to prove the following existence theorem for QFDE \eqref{fie.}.
\begin{thm}\label{c}
Assume that  the assumption $(\mathcal{H}_0)$) through $(\mathcal{H}_7)$ hold. Furthermore, if
\begin{equation}\label{equa1}
\left\{
\begin{array}{ll}
   \displaystyle L\left(\left|\frac{x_0- f_1(0,x_0)}{f_2(0,y_0)}\right| + \|h\|_{L^1}\right)\leq K & \hbox{} \\\\
   M_1\|\gamma\|+K\leq 1-k.& \hbox{}
   \end{array}
   \right.
   \end{equation}
Then the system of the functional differential  equations \eqref{fie.} has, at least, one solution in
$C(J,\R)\times C(J,\R)$.$\hfill\diamondsuit$
\end{thm}
\begin{pf}
Observe that the above problem \eqref{fie.} may be written in the following form
\begin{equation*}
\left\{
  \begin{array}{ll}
    x(t)=(Ax)(t)+(By)(t)\cdot(B'y)(t) & \hbox{} \\
    y(t)=(Cx)(t)+(Dy)(t) & \hbox{}.
  \end{array}
\right.
\end{equation*}
where $A$, $B$, $C$, $D$ and $B'$ on $C(J,\R)$ defined by:
\begin{equation*}
\left\{
  \begin{array}{ll}
    (Ax)(t) = f_1(t,x(t)) & \hbox{} \\  \\
    (By)(t) = f_2(t,y(t)) & \hbox{}\\\\
    (Cx)(t) = \displaystyle\frac{1}{1-b(t)|x(t)|}-p\left(t,\frac{1}{1-b(t)|x(t)|}\right) \\ \\
    (Dy)(t) = p(t,y(t)) & \hbox{} \\ \\
    (B'y)(t) = \displaystyle c e^{-\lambda t}+ e^{-\lambda t}\int_{0}^{t}e^{\lambda s}g(s,y(s))ds.
  \end{array}
\right.
\end{equation*}
From the above assumptions and the continuity of the integral, it follows that the operators $B,B':E\longrightarrow\mathcal{K}.$
In order to apply Theorem \ref{th1}, we have to verify the following steps.\\
$\mathbf{Step~1}$: $A$, $B$, $C$, $D$ and $B'$ are nondecreasing on $E$.\\
Let $x,y \in E$ be such that $x\geq y$. Then $x(t)\geq y(t)$ for all $t\in J$. Then by assumption $(\mathcal{H}_3)$, we obtain
\begin{equation*}
\begin{array}{llll}
Ax(t) &=&f_1(t,x(t)\\
      &\geq& f_1(t,y(t)\\
      &=& Ay(t),
\end{array}
\end{equation*}
for all $t \in J$. This shows that the operator that the operator $A$ is nondecreasing on $E$. Similarly, by assumption $(\mathcal{H}_3)$, we get
\begin{equation*}
\begin{array}{llll}
Bx(t) &=&f_2(t,x(t)\\
      &\geq& f_2(t,y(t))\\
      &=& By(t),
\end{array}
\end{equation*}
for all $t \in J$. This shows that the operator  $B$ is also nondecreasing on $E$. Furthermore, by assumption $(\mathcal{H}_5)$, we get $C$ is nondecreasing operator on $E$. Indeed, let $x,y\in E$ such that $x(t)\geq y(t)$
\begin{equation*}
\begin{array}{llllllll}
Cx(t)&=&\displaystyle\frac{1}{1-b(t)|x(t)|}-p\left(t,\frac{1}{1-b(t)|x(t)|}\right)\\
      &\geq&\displaystyle\frac{1}{1-b(t)|y(t)|}-p\left(t,\frac{1}{1-b(t)|y(t)|}\right)\\
      &=&Cy(t),
            \end{array}
\end{equation*}
for all $t \in J$. This shows that the operator $C$ is nondecreasing on $E$. Again, by
assumption $(\mathcal{H}_4)$ , we obtain
\begin{equation*}
\begin{array}{llll}
Dx(t) &=&p(t,x(t)\\
      &\geq& p(t,y(t))\\
      &=& Dy(t),
\end{array}
\end{equation*}
for all $t \in J$. This shows that the operator  $D$ is nondecreasing on $E$. Finally, by
assumption $(\mathcal{H}_5)$ , we get
\begin{equation*}
\begin{array}{llll}
B'x(t) &=&\displaystyle c e^{-\lambda t}+ e^{-\lambda t}\int_{0}^{t}e^{-\lambda s}g(s,x(s))ds\\
      &\geq&\displaystyle ce^{-\lambda t}+ e^{-\lambda t}\int_{0}^{t}e^{-\lambda s}g(s,y(s))ds\\
      &=& B'y(t),
\end{array}
\end{equation*}
for all $t \in J$. This shows that the operator  $B'$ is nondecreasing on $E$.\\
$\mathbf{Step~2}$: $A$, $B$ and $C$ are partially bounded and partially $\mathcal{D}$-Lipschitzians on $E$.\\
Let $x \in E$ be arbitrary. Without loss of generality we assume that $x \geq 0$. Then by assumptions $(\mathcal{H}_1)$ and $(\mathcal{H}_2)$, we have
\begin{equation*}
\begin{array}{llll}
|Ax(t)| &=&  |f_1(t,x(t))|\\
        &\leq& |f_1(t,x(t))-f_1(t,0)|+ |f_1(t,0)|\\
        &\leq&\displaystyle \frac{L\|x\|}{K+\|x\|}+ F_0\\
        &\leq& L+F_0,
        \end{array}
\end{equation*}
for all $t \in J$. Taking the supremum over $t$ in the above inequality, we obtain
$$\|Ax\|\leq L+F_0$$
for all $x \in E$. So, $A$ is bounded. This further implies that $A$ is partially bounded on $E$.\\
Next, let $x,y \in E$ be such that $x \geq y$. Then, we have
\begin{equation*}
\begin{array}{llll}
|Ax(t)-Ay(t)| &=& |f_1(t,x(t))-f_1(t,y(t))|\\
              &\leq& \displaystyle \frac{L\|x-y\|}{K+\|x-y}\|\\
              &=&\psi_A(\|x-y\|),
\end{array}
\end{equation*}
for all $t\in J$, where $\displaystyle\psi_A(r)= \frac{Lr}{K+r}$. Taking the supremum over $t$, we obtain
$$\|Ax-Ay\|\leq \psi_A(\|x-y\|),$$
for all $x,y \in E$ with $x \geq y$. Hence, $A$ is a partial nonlinear $\mathcal{D}$-Lipschitzian on $E$ with a $\mathcal{D}$-function $\psi_A$.
By using the same argument, we conclude that $B$ is partially bounded and partially $\mathcal{D}$-Lipschitzian
on $E$, where $\displaystyle\psi_B(r)=\frac{Lr}{K+r}$.
Also, we shall show that $C$ is partially bounded and partially $\mathcal{D}$-Lipschitzian. Indeed, for all $t\in J$, we get
\begin{eqnarray*}
\displaystyle
|Cx(t)| &=& \left|\frac{1}{1-b(t)|x(t)|}- p\left(t,\frac{1}{1-b(t)|x(t)|}\right)\right|\\
         &\leq& 1 +  \left|p\left(t,\frac{1}{1-b(t)|x(t)|}\right)\right|\\
         &\leq& 1+k
\end{eqnarray*}
This means that the operator $C$ is partially bounded. Moreover, let $x,y \in E$ such that $x\geq y$. Then we get
\begin{equation*}
\begin{array}{llllllll}
|Cx(t)-C(y)(t)|&\leq&\displaystyle\left|\frac{1}{1-b(t)|x(t)|}- p\left(t,\frac{1}{1-b(t)|x(t)|}\right)-
\frac{1}{1-b(t)|y(t)|}-p\left(t,\frac{1}{1-b(t)|y(t)|}\right)\right|\\ \\
&\leq&\displaystyle\left|\frac{1}{1-b(t)|x(t)|}-\frac{1}{1-b(t)|y(t)|}\right|+
\left|p\left(t,\frac{1}{1-b(t)|x(t)|}\right)-p\left(t,\frac{1}{1-b(t)|y(t)|}\right)\right|\\ \\
            &\leq& (1+k)\displaystyle\left|\frac{1}{1-b(t)|x(t)|}-\frac{1}{1-b(t)|y(t)|}\right|\\\\
            &\leq&(1+k)\|b\||x(t)-y(t)|.
\end{array}
\end{equation*}
Taking the supremum over $t$, we obtain that $C$ is partially nonlinear $\mathcal{D}$-Lipshitzian with $\mathcal{D}$-function $\psi_C(r)= (1+k)\|b\|r$\\
$\mathbf{Step~3}$ : $(I-D)^{-1}$ is partially nonlinear $\mathcal{D}$-Lipshitzian and $(I-D)^{-1}C$ is nondecreasing.\\
Since $D$ is contraction then $(I-D)^{-1}$ exist and is a  contraction with constant $\displaystyle\frac{1}{1-k}$. Consequently, $(I-D)^{-1}$ is partially $\mathcal{D}$-lipshitzian with $\mathcal{D}$-function $\displaystyle\psi_\phi(r)=\frac{1}{1-k}r$\\
Now, we shaw that $(I-D)^{-1}C$  is nondecreasing. Since $\displaystyle(I-D)^{-1}Cx=\frac{1}{1-b|x|}$, for all $x\in E$. Then for all $x,y\in E$ such that $x\leq y$ we have
\begin{eqnarray*}
\displaystyle(I-D)^{-1}Cx=\frac{1}{1-b|x|}
                         \leq\frac{1}{1-b|y|}
                         =(I-D)^{-1}Cy.
\end{eqnarray*}
$\mathbf{Step~4}$ : $B'$ est partially continuous and $C$ is compact.\\
Let $\{x_n\}_{n\in \N}$  be a sequence in a chain $C$ such that $x_n \to x$ as $n \to
\infty$. Then by the dominated convergence theorem for integration, we obtain
\begin{equation*}
\begin{array}{llll}
\displaystyle\lim_{n\longrightarrow \infty}Bx_n(t) &=& ce^{-\lambda t} + e^{-\lambda t}\displaystyle\int_{0}^{t}
                                          \lim_{n\longrightarrow \infty}e^{\lambda s}g(s,x_n(s))ds \\
                                      &=& ce^{-\lambda t} + e^{-\lambda t}\displaystyle\int_{0}^{t} e^{\lambda
                                      s}g(s,x(s))ds \\
                                      &=& Bx(t),
\end{array}
\end{equation*}
for all $t \in J$. This shows that $\{Bx_n\}$ converges to $Bx$ pointwise on $J$.\\
Now we show that $\{Bx_n\}_{n\in \N}$ is an equicontinuous sequence of functions in $E$. Let $t_1, t_2 \in J$
with $t_1 > t_2 > 0$. Then we have
\begin{equation*}
\begin{array}{lllll}
  |Bx_n(t_2)-Bx_n(t_1)| &\leq& \displaystyle\left|ce^{-\lambda t_2}-ce^{-\lambda t_1}\right|+\left|e^{-\lambda
  t_2}\int_0^{t_2}e^{\lambda s}g(s,x_n(s))ds\right. \\
  &-&\left.e^{-\lambda t_1}\int_0^{t_1}e^{\lambda s}g(s,x_n(s))ds\right|  \\
   &\leq& \displaystyle|ce^{-\lambda t_2}-ce^{-\lambda t_1}|+\left|e^{-\lambda t_2}\int_0^{t_2}e^{\lambda
   s}g(s,x_n(s))\right.ds\\
   &-&\displaystyle \left.e^{-\lambda t_2}\int_0^{t_1}e^{\lambda s}g(s,x_n(s))ds\right.
   \displaystyle\left.+e^{-\lambda t_2}\int_0^{t_1}e^{\lambda s}g(s,x_n(s))ds\right.\\
   &-&\left.\displaystyle e^{-\lambda t_1} \int_0^{t_1}e^{\lambda s}g(s,x_n(s))ds\right|\\
   &\leq& \displaystyle |ce^{-\lambda t_2}-ce^{-\lambda t_1}|+\left|e^{-\lambda t_2}\int_{t_1}^{t_2}e^{\lambda
   s}g(s,x_n(s))ds\right|\\
    &+&\displaystyle\left|(e^{-\lambda t_2}-e^{-\lambda t_1}) \int_0^{t_1}e^{\lambda s}g(s,x_n(s))ds
    \right|.
\end{array}
\end{equation*}
This implies that $\displaystyle \lim_{t_1\to t_2}|Bx_n(t_2)-Bx_n(t_1)|=0$ and consequently $B$ is partially continuous.\\
Next, we show that $C$ is compact. To do this, for any $x\in J$
\begin{eqnarray*}
\displaystyle
|Cx(t)| &=& \left|\frac{1}{1+b(t)|x(t)|}- p\left(t,\frac{1}{1+b(t)|x(t)|}\right)\right|\\
         &\leq& 1 +  \left|p\left(t,\frac{1}{1+b(t)|x(t)|}\right)\right|\\
         &\leq& 1+k
\end{eqnarray*}
this show that $C$ is bounded. Now we show that $C$ is equicontinuous on $E$.
For each $t_1, t_2 \in J$ such that $t_2>t_1$, we get
\begin{eqnarray*}
\displaystyle
|Cx(t_2) - Cx(t_1)| &=& \left|\frac{1}{1+b(t)|x(t_2)|}- p\left(t_2,\frac{1}{1+b(t_2)|x(t_2)|}\right)\right.\\
&&-\left.\frac{1}{1+b(t_1)|x(t_1)|}-p\left(t_1,\frac{1}{1+b(t_1)|x(t_1)|}\right)\right|\\
&\leq &\left|\frac{1}{1+b(t)|x(t_2)|} -\frac{1}{1+b(t_1)|x(t_1)|}\right| \\
&& + \left|p\left(t_2,\frac{1}{1+b(t_2)|x(t_2)|}\right)
-p\left(t_1,\frac{1}{1+b(t_1)|x(t_1)|}\right) \right|\\
&\leq&\|b\|_\infty\left|\frac{x(t_2)-x(t_1)}{(1+b(t_2)|x(t_2)|)(1+b(t_1)|x(t_1)|)}\right|\\
&&+k\|b\|_\infty \left|\frac{x(t_2)-x(t_1)}{(1+b(t_2)|x(t_2)|)(1+b(t_1)|x(t_1)|)} \right|\\
&\leq& (1+k)\|b\|_\infty |x(t_2) - x(t_1)|.
\end{eqnarray*}
This implies that $$\lim_{t_2\to t_1}|Cx(t_2)-Cx(t_1)| = 0$$
This means that $C$ is equicontinuous on J. Then by the Arzela-Ascoli's theorem \cite{curtain} the closure of $C(E)$ is relatively compact, consequently $C$ is compact operator.\\
$\mathbf{Step~5}$: $u$ satisfies the operator inequality $u\leq Au+B(I-D)^{-1}Cu\cdot B'(I-D)^{-1}Cu$.\\
By assumption $(\mathcal{H}_6)$, we have
\begin{equation*}
\displaystyle
\left\{
\begin{array}{llll}
\displaystyle\left(\frac{u(t)-f_1(t,u(t))}{f_2(t,v(t))}\right)'+ \lambda
\left(\frac{u(t)-f_1(t,u(t))}{f_2(t,v(t))}\right)
\leq g(t,v(t))\\
\displaystyle v(t) \leq \displaystyle\frac{1}{1-b(t)|u(t)|}-p\left(t,\frac{1}{1-b(t)|u(t)|}\right)+p(t,v(t))\\ \\
\big(u(0),v(0)\big)\leq(u_0,v_0) \in \R^2
\end{array}\right.
\end{equation*}
for all $t\in J$.

Then, by Remark \ref{lem1} we get
\begin{equation*}
\displaystyle
\left\{
\begin{array}{llll}
u(t)&\leq& f_1(t,u(t))+f_2(t,v(t))\left(ce^{-\lambda t}+e^{-\lambda t}\displaystyle\int_{0}^{t}e^{\lambda s}g(s,v(s))ds\right)\\
\displaystyle v(t) &\leq& \displaystyle\frac{1}{1-b(t)|u(t)|}-p\left(t,\frac{1}{1-b(t)|u(t)|}\right)+p(t,v(t))\\
\end{array}
\right.
\end{equation*}
Hence, from definitions of the operators $A,B,C,D$ and $B'$ it follows that
$$u(t)\leq Au(t) + B(I-D)^{-1}Cu(t)\cdot B'(I-D)^{-1}Cu(t),$$
for all $t\in J$. Taking the suprumum we have
$$u\leq Au + B(I-D)^{-1}Cu\cdot B'(I-D)^{-1}Cu.$$

 Finally  By using the assumption $(\mathcal{H}_5)$, we have
\begin{eqnarray*}
M & = & \|T'(E)\|\\
  & = & \displaystyle \left\|ce^{-\lambda t} + e^{-\lambda t}\int_{0}^{t}e^{\lambda t}g(s,x(s))ds\right\|\\
  & \leq & |c| + \|h\|_{L^1}.
\end{eqnarray*}
From equation \eqref{equa1}, we infer that $M\psi_B\circ\psi_\phi\circ\psi_C(r)+\psi_A(r)<r$.\\
Thus, the operators  $A,B,C,D$ and $B'$ satisfy all the requirement of Theorem \ref{th1} and so the QFDE
\eqref{fie.} has a positive solution in $C(J\times\R)\times C(J\times\R)$.
\end{pf}
\begin{rem}
The conclusion of Theorem \ref{c} also remains true if we replace the assumption $(\mathcal{H}_7)$ with the following one:\\
$(\mathcal{H}'_7)$ The QFDE \eqref{fie.} has a upper solution in $C(J,\R)\times C(J,\R)$
\end{rem}
 The proof under this new assumption is similar to Theorem \ref{c} and the conclusion again follows by an application of Theorem \ref{th1}.

\end{document}